\documentclass[conference]{IEEEtran}
\pdfoutput=1
\usepackage{cite}
\usepackage{graphicx}
\usepackage{amsmath}
\usepackage{array}

\usepackage[square,sort,comma,numbers]{natbib}
\usepackage{subcaption}
\begin{document}
\title{Higher Order Mandelbrot Sets and their Varying Shapes}
\author{\IEEEauthorblockN{Arshdeep Singh Pareek}}

\maketitle

\begin{abstract}
Mandelbrot set arose from the pioneering work of French mathematician Gaston Julia in the field of complex dynamics at the beginning of the 20th century. French-American mathematician Benoit Mandelbrot used computers to calculate iterations of complex polynomials of second order and displayed intricate images of fractal geometry. While studying fundamental properties of the Mandelbrot set, little attention has been paid to study the relationship between the degree of the generating polynomials and the shape of the main body of the Mandelbrot set. This paper extends the work from generating polynomials of second degree to polynomials of higher degrees using basic principles of complex numbers and calculus and shows that the number of primary lobes in the Mandelbrot set corresponding to a polynomial of nth degree is n-1 and as n tends to infinity the shape of the primary lobe tends to a circle of unit radius. It is also shown that the arguments of the points where secondary lobes are attached to the primary lobes are given by the roots of unity. These results provide an easy way to predict shapes of the main lobes of Mandelbrot sets and the locations of the points where secondary lobes are attached to the main lobe and may be helpful in understanding features of the Mandelbrot set.
\end{abstract}

\section{Introduction}

 Mandelbrot set is defined based on the idea of iterations of complex polynomial functions of second degree. A complex function is repeatedly applied to some initial arguments and  the behaviour of the resulting values is observed. For example, for a function of the form $f(z): z^2 + c$ defined in the complex plane where $z = x + iy$, and $c$ is a complex number, different values of $c$ would produce widely different sequences of $z$. The sequences can be categorized as either bounded, where the values of the sequence stay within an arbitrary bound, or unbounded, where the values can be made larger than any arbitrary number as the number of iterations increase.

 This idea of the set of complex numbers $c$ in an iterative function of the form $f(z): z^n + c$, where the resulting sequences are bounded, can be extended to cases where the value of $n$ is greater than 2. Such sets are termed Multibrot sets which show interesting behaviour as the value of $n$ increases.

\subsection{Orbits}
To better understand the behaviour of sequences resulting from iterative complex valued functions, we can plot them on the two dimensional complex plane. Let us first consider the simple real function $f(z): z^2$ which is a special case of $f(z): z^n + c$ where $n$ is 2 and $c$ equals 0. 

\begin{figure*}
    
    \begin{subfigure}[b]{1\columnwidth}
        \includegraphics[width=0.7\textwidth]{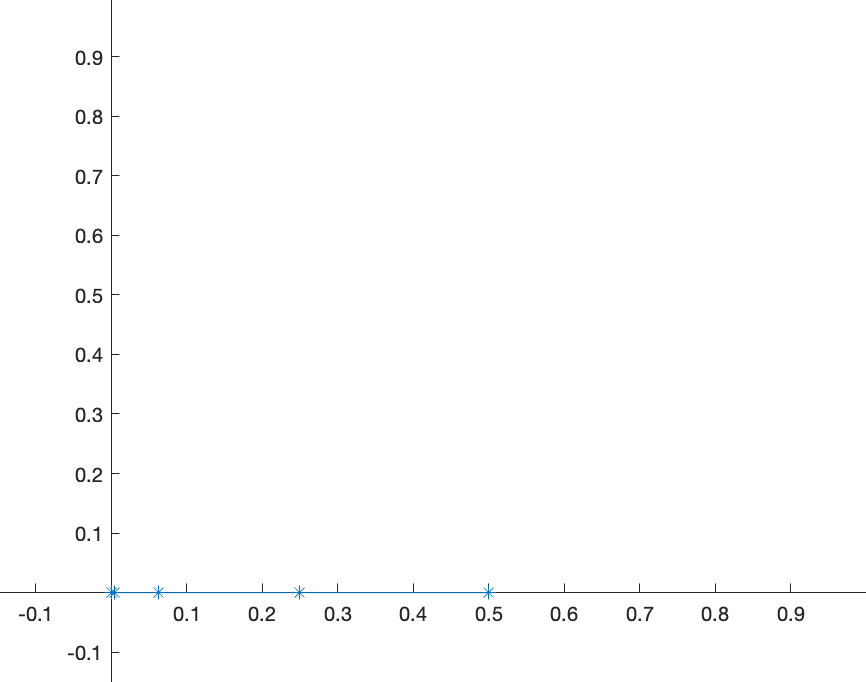}
        \caption{Orbit for c = 0, z = 0.5.}
        \label{orbit_1}
    \end{subfigure}
    \hfill
    \begin{subfigure}[b]{1\columnwidth}
        \includegraphics[width=0.7\textwidth]{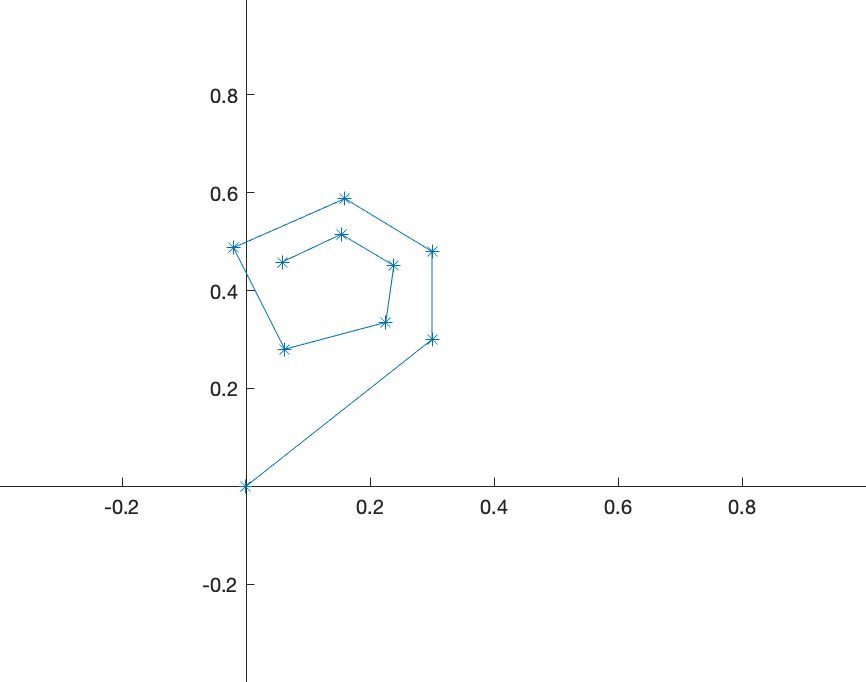}
        \caption{Orbit for c = 0.3 + 0.3i, z = 0.}
        \label{orbit_2}
    \end{subfigure}

    \begin{subfigure}[b]{1\columnwidth}
        \includegraphics[width=0.7\textwidth]{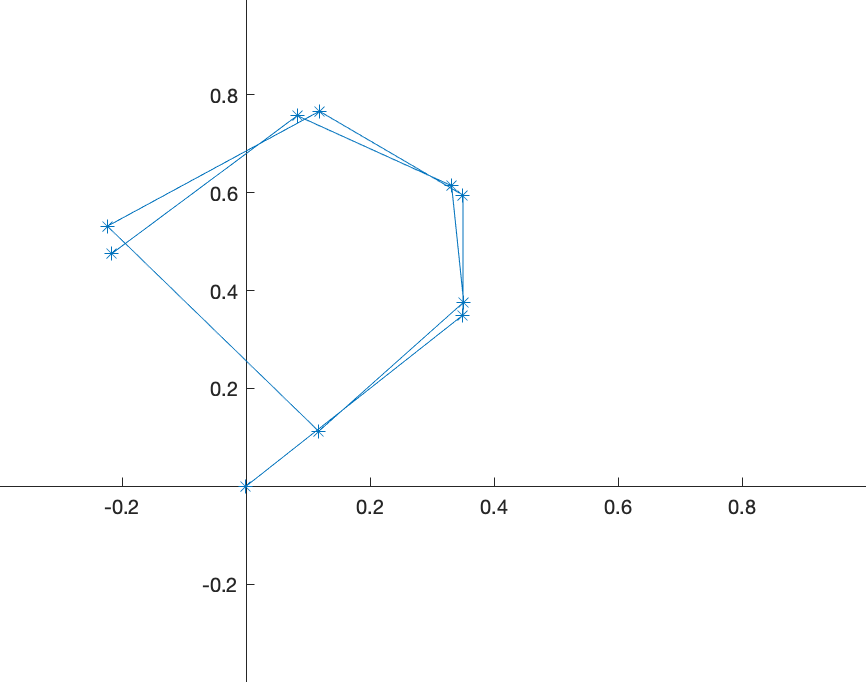}
        \caption{Orbit for c = 0.35 + 0.35i, z = 0.}
        \label{orbit_3}
    \end{subfigure}
    \hfill
    \begin{subfigure}[b]{1\columnwidth}
        \includegraphics[width=0.7\textwidth]{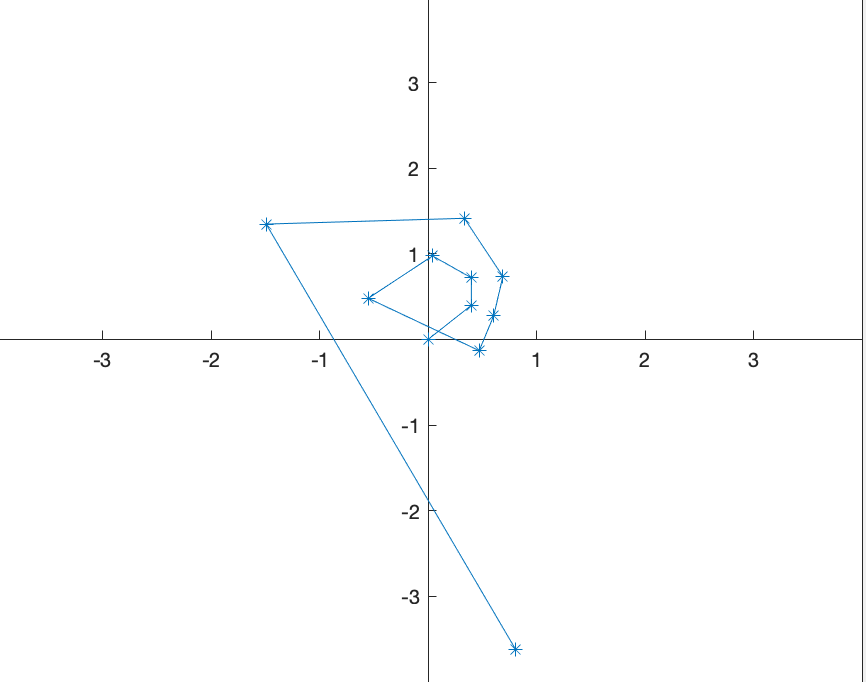}
        \caption{Orbit for c = 0.4 + 0.4i, z = 0.}
        \label{orbit_4}
    \end{subfigure}
    \caption{Orbits}
    \label{fig:Orbits}
\end{figure*}

When $z = 2$ is taken as the initial value, the repeated iterations of the function results in the sequence 2, 4, 16, 256, and so on. In case the iterations are started with $z = 0.5$ (Figure \ref{orbit_1}) as the initial point, the successive iterations are 0.25, 0.125, 0.0625, and so on. It is seen that if z = 2 (or any value greater than 1) is taken as the starting point, each resulting successive point in the sequence increases in its value than the preceding point. In mathematical terms such sequences are called unbounded and the orbits of such points escape to infinity. On the contrary if the iterations are started with initial value $z = 0.5$
(or with any value less than 1), the successive points generated by iterations grow smaller and smaller in value and
approach zero. The iterations of the initial point z=1 remain bounded.

  The term ``orbits" becomes more apparent when the sequences are observed for the function $f(z): z^2 + c$ as the value of $c$ is made non-zero. For example, when $c = 0.3 + 0.3i$ the graph of the function values resembles an orbit around an attractor point. This can be observed in Figure \ref{orbit_2}. While the orbit in Figure \ref{orbit_2} tends towards a single point and is bounded. Figure \ref{orbit_4}, on the other hand, shows an orbit which is unbounded.

\subsection{Fixed Points}
 
A fixed point of a function is a point, say $w$, which yields the same output. In other words, $f(w) = w$.\cite{peitgen1992} Fixed points can be further divided into three categories, attractor, repellor and neutral. For
example, for function $f(z) = z^2$, $z = 0$ is an attractor as any point z in the vicinity of 0 will be ‘pulled’ towards 0. For the same function, $z = 1$ is a repellor as iterations of function $f$ for any value of $z$ in the
vicinity of 1 will be ‘repelled’ either towards 0 (if $|z| < 1$) or towards infinity (if $|z| > 1$).

For a fixed point $w$, if $|f’(w)| < 1$, the fixed point is attractive. If $|f’(w)| > 1$, the fixed point w is repelling. In case $|f’(w)|=1$, the fixed point is classified as neutral.

\section{Mandelbrot set}

The Mandelbrot set (M-set) of a given function $f$ is defined as the set of all $c$ for a given value of $n$
for which the orbits of initial points $z_0 = 0$ (critical orbits) are bounded. In other words M-set can be
defined as $M = {c \in C : f_c^k(0)\ remains\ finite\ as\ k\xrightarrow[]{} \infty}$.\cite{peitgen1992}

With the above background knowledge, Mandelbrot sets of polynomials of different degrees can be examined. Basic properties of the Mandelbrot set can be derived applying its definition and application of simple mathematical tools.

\subsection{Mandelbrot set for $f(z) = z^2 + c$}
The first and perhaps the easiest to tackle is the Mandelbrot set generated by a quadratic polynomial function $f(z) = z^2 + c$.

%\begin{figure}[h]
%    \centering
%    \includegraphics[width=0.4\textwidth]{MandelbrotSet.png}
%    \caption{The Mandelbrot Set}
%    \label{fig:MandelbrotSet}
%\end{figure}
    
It must be kept in mind that the Mandelbrot Set is not a trace of a continuous function, but a collection of discrete points, ‘c’, in a complex plane which generates sets of bounded orbits from the function, $f(z) = z^2 + c$.

One of the properties which can be easily derived is the fact that the Mandelbrot set is confined to a limited region. For the quadratic polynomial, $f(z) = z^2 + c$, let us examine the orbit of $z_0 = 0$. Let $z_n$ be the nth iterate of $z_0$ whose absolute value is not less than the modulus of $c$, and greater than 2. For $c$ to belong to the Mandelbrot set, the orbit of $z_0$ must be bounded. If there exists an $n$ for which $|z_n| > 2$. Then there exists a small but positive number $\epsilon > 0$ such that $|z_n| = 2 + \epsilon$. For the starting point $z_0 = 0$ of the iteration, it is assumed that $|z_n| > |c|$ and $|z_n| > 2$ for some $n$.

Triangle inequality for complex numbers imply,

\begin{equation}
\begin{aligned}
|z_n^2| = |z_n^2 + c - c| \leq |z_n^2 + c| + |c|\nonumber
\end{aligned}
\end{equation}
\begin{equation}
\begin{aligned}
\implies |z_n^2 + c| \geq |z_n^2| - |c| \geq |z_n|^2 - |z_n|\nonumber
\end{aligned}
\end{equation}
\begin{equation}
\begin{aligned}
= \left(|z_n|\ -\ 1\right)|z_n| = \left(2 + \epsilon\ -\ 1\right) |z_n| = \left(1+\epsilon\right)|z_n|
\end{aligned}
\end{equation}

Or, $z_{n+1} > z_n$

Thus, if the function f is iterated, each iteration will increase the absolute value by $\left(1 + \epsilon\right)$ and after $k$ iterations the value will be $\left(1 + \epsilon\right)^k |z_n|$ which will tend to infinity as $k$ tends to infinity. In other words the sequence of iterations of the initial point $z_0 = 0$ is not bounded. It follows from this proposition that for points $z$ such that $|z| = max\left(|c|, 2\right)$ the orbits will fly to infinity and accordingly all such point will not be part of the M-set.\cite{peitgen1992}

\subsection{Shape of the main body of the Mandelbrot set for a quadratic function}

The shape of the central part of the Mandelbrot set of the quadratic polynomial can be plotted using parametric equations. The following is the analytical method to evolve the parametric equations.
For $f(z) = z^2 + c$; let $z$ be a fixed point. Then by definition of a fixed point,

\begin{equation}
\begin{aligned}
f(z) = z\nonumber
\end{aligned}
\end{equation}
\begin{equation}
\begin{aligned}
\implies z^2 + c = z\nonumber
\end{aligned}
\end{equation}
\begin{equation}
\begin{aligned}
\implies z^2\ -\ z + c = 0
\end{aligned}
\label{eq.1}
\end{equation}

Since the Mandelbrot set will include all the attractors (i.e. $|f'(z)| < 1$ as mentioned in the definition of fixed points) and exclude repellors (i.e. $|f'(z)| > 1$) while determining the shape of the main part of the Mandelbrot sets, one uses the criterion $|f'(z)| = 1$ as it represents the transition between the attracting points and repelling points which is the boundary of the Mandelbrot set.

Thus,
\begin{equation}
\begin{aligned}
f'(z) = 1\nonumber
\end{aligned}
\end{equation}
\begin{equation}
\begin{aligned}
\implies 2z = 1\nonumber
\end{aligned}
\end{equation}
\begin{equation}
\begin{aligned}
\implies 2z = cis\phi
\end{aligned}
\end{equation}

Substituting $z = cis\phi/2$ in Equation \ref{eq.1},

\begin{equation}
\begin{aligned}
\frac{cis2\phi}{4} - \frac{cis\phi}{2} + c = 0\nonumber
\end{aligned}
\end{equation}
\begin{equation}
\begin{aligned}
\implies c = \frac{cis\phi}{2} - \frac{cis2\phi}{4}
\end{aligned}
\end{equation}

Representing $c$ as $x + iy$ and separating real and imaginary:

\begin{equation}
\begin{aligned}
x + iy = \frac{cos\phi}{2} - \frac{cos2\phi}{4} + i\left(\frac{sin\phi}{2} - \frac{sin2\phi}{4}\right)\nonumber
\end{aligned}
\end{equation}
\begin{equation}
\begin{aligned}
x = \frac{cos\phi}{2} - \frac{cos2\phi}{4}\nonumber
\end{aligned}
\end{equation}
\begin{equation}
\begin{aligned}
y = \frac{sin\phi}{2} - \frac{sin2\phi}{4}
\end{aligned}
\end{equation}

This is a set of parametric equations in $\phi$ which can take values from 0 radians onwards.\cite{peitgen1992} A plot of this set of parametric equations is clearly a cardioid, as can be seen in Figure \ref{fig:card}. It can be easily seen that when $\phi=0$, $x = \frac{1}{4}$ and when $\phi=\pi$, $x = -\frac{3}{4}$. In other words, the $x$ intercepts of the cardioid are at $x = -\frac{3}{4}$. and $x = \frac{1}{4}$.

\subsection{Shape of the main body of the Mandelbrot set for a cubic function}
The same method will now be used to determine the shape of the main body of the Mandelbrot set of the cubic polynomial, $f(z) = z^3+c$;

Let z be a fixed point. Then by definition of the fixed point,

\begin{equation}
\begin{aligned}
f(z) = z\nonumber
\end{aligned}
\end{equation}
\begin{equation}
\begin{aligned}
\implies z^3 + c = z\nonumber
\end{aligned}
\end{equation}
\begin{equation}
\begin{aligned}
\implies z^3 - z + c = 0\nonumber
\end{aligned}
\end{equation}
\begin{equation}
\begin{aligned}
\implies c = z - z^3
\end{aligned}
\label{eq(c)}
\end{equation}

For a fixed point the derivative of the function should be equal to 1.

Thus,

\begin{equation}
\begin{aligned}
f'(z) = 1\nonumber
\end{aligned}
\end{equation}
\begin{equation}
\begin{aligned}
\implies 3z^2 = 1\nonumber
\end{aligned}
\end{equation}
\begin{equation}
\begin{aligned}
\implies 3z^2 = cis\phi
\end{aligned}
\end{equation}

Substituting $z = \frac{cis\frac{\phi}{2}}{\sqrt{3}}$ in Equation \ref{eq(c)}, and representing $c$ as $x+iy$ and separating real and imaginary parts:

\begin{equation}
\begin{aligned}
x + iy = \frac{cos\frac{\phi}{2}}{\sqrt{3}} - \frac{cos\frac{3\phi}{2}}{3\sqrt{3}} + i\left(\frac{sin\frac{\phi}{2}}{\sqrt{3}} - \frac{sin\frac{3\phi}{2}}{3\sqrt{3}}\right)\nonumber
\end{aligned}
\end{equation}
\begin{equation}
\begin{aligned}
x = \frac{cos\frac{\phi}{2}}{\sqrt{3}} - \frac{cos\frac{3\phi}{2}}{3\sqrt{3}}\nonumber
\end{aligned}
\end{equation}
\begin{equation}
\begin{aligned}
y = \frac{sin\frac{\phi}{2}}{\sqrt{3}} - \frac{sin\frac{3\phi}{2}}{3\sqrt{3}}
\end{aligned}
\end{equation}

This is a set of parametric equations in $\phi$ which can take values from 0 radians onwards. A plot of this set of parametric equations is clearly a double-cardioid as can be seen in Figure \ref{fig:card}.

\subsection{Shape of the main body of the Mandelbrot set for a polynomial of nth degree}

The shape of the central part of the Mandelbrot set of the polynomial of nth degree will now be determined, $f(z) = z^n+c$; let $z$ be a fixed point. Then by definition of the fixed point,

\begin{equation}
\begin{aligned}
f(z) = z\nonumber
\end{aligned}
\end{equation}
\begin{equation}
\begin{aligned}
z^n+c = z\nonumber
\end{aligned}
\end{equation}
\begin{equation}
\begin{aligned}
z^n - z + c = 0\nonumber
\end{aligned}
\end{equation}
\begin{equation}
\begin{aligned}
c = z - z^n
\end{aligned}
\label{eqcn}
\end{equation}

For a fixed point the derivative of the function should be equal to 1.

Thus,

\begin{equation}
\begin{aligned}
f'(z) = 1\nonumber
\end{aligned}
\end{equation}
\begin{equation}
\begin{aligned}
nz^{n-1} = 1\nonumber
\end{aligned}
\end{equation}
\begin{equation}
\begin{aligned}
nz^{n-1} = cis\phi\nonumber
\end{aligned}
\end{equation}

Substituting $z =\left(\frac{cis\phi}{n}\right)^{\frac{1}{n-1}}$ in equation \ref{eqcn},

\begin{equation}
\begin{aligned}
c = \left(\frac{cis\phi}{n}\right)^{\frac{1}{n-1}} - \left(\frac{cis\phi}{n}\right)^\frac{n}{n-1}
\end{aligned}
\end{equation}

Representing c as x+iy and separating real and imaginary parts,

\begin{equation}
\begin{aligned}
x + iy = \frac{cos\frac{\phi}{n-1}}{n^{\frac{1}{n-1}}} - \frac{cos\frac{n\phi}{n-1}}{n^{\frac{n}{n-1}}} + i\left(\frac{sin\frac{\phi}{n-1}}{n^{\frac{1}{n-1}}} - \frac{sin\frac{n\phi}{n-1}}{n^{\frac{n}{n-1}}}\right)\nonumber
\end{aligned}
\end{equation}
\begin{equation}
\begin{aligned}
x = \frac{cos\frac{\phi}{n-1}}{n^{\frac{1}{n-1}}} - \frac{cos\frac{n\phi}{n-1}}{n^{\frac{n}{n-1}}}\nonumber
\end{aligned}
\end{equation}
\begin{equation}
\begin{aligned}
y = \frac{sin\frac{\phi}{n-1}}{n^{\frac{1}{n-1}}} - \frac{sin\frac{n\phi}{n-1}}{n^{\frac{n}{n-1}}}
\end{aligned}
\end{equation}

This is set of parametric equations in $\phi$ which can take values from $0$ to $2\pi$ radians.

The above parametric equations are functions of sine and cosine of angles $\frac{\phi}{n-1}$ which have periods of $2(n-1)\pi$. In other words, for n=2, the period is $2\pi$, for n=3 the period is $4\pi$, for n=4 the period is $6\pi$ and so on.

\begin{figure}[h!]
  \centering
  \includegraphics[width=.8\linewidth]{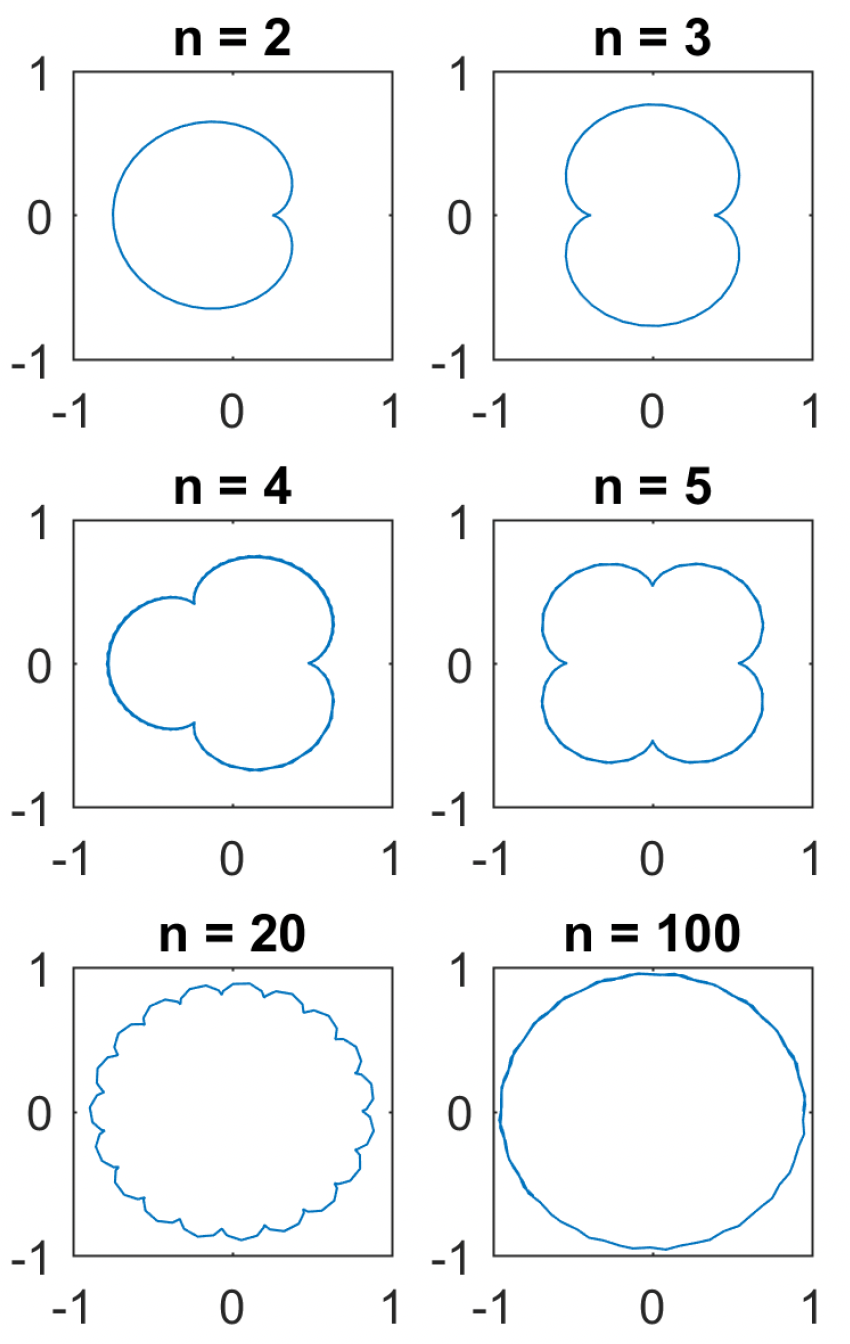}
  \caption{Cardioids.}
  \label{fig:card}
\end{figure}

The computer program MATLAB has been used to plot the graphs of these parametric equations for various values of n. Graphs for n=2, 3, 4, 5, 20 and 100 are given in Figure \ref{fig:card}.

\section{Mandelbrot set - Convergence to a unit circle for large values of n}

As can be seen in the figure above, the Mandelbrot sets of increasing power seem to tend towards a circle; specifically a unit circle. Hence, in order to further examine the shapes of the plots of Mandelbrot sets for large values of $n$, it would be useful to find the value of $x^2 + y^2$, which is the square of the distance of the point $(x,y)$ from the origin, where $x$ and $y$ are given by the parametric equations above, and examine how does this behave. Squaring and adding the values of x and y from the above parametric equations and using trigonometric identities:

\begin{equation}
\begin{aligned}
sin^2\theta + cos^2\theta = 1\nonumber
\end{aligned}
\end{equation}

\begin{equation}
\begin{aligned}
cosAcosB + sinAsinB = cos\left(A - B\right)\nonumber
\end{aligned}
\end{equation}
\vfill

\pagebreak
One gets,
\begin{equation}
\begin{aligned}
x^2+y^2 = \left(\frac{cos\frac{\phi}{n-1}}{n^{\frac{1}{n-1}}} - \frac{cos\frac{n\phi}{n-1}}{n^{\frac{n}{n-1}}}\right)^2 + \left(\frac{sin\frac{\phi}{n-1}}{n^{\frac{1}{n-1}}} - \frac{sin\frac{n\phi}{n-1}}{n^{\frac{n}{n-1}}}\right)^2\nonumber
\end{aligned}
\end{equation}
\begin{equation}
\begin{aligned}
x^2 + y^2 = \left(\frac{cos^2\frac{\phi}{n-1} + sin^2\frac{\phi}{n-1}}{\left(n^{\frac{1}{n-1}}\right)^2}\right) + \left(\frac{cos^2\frac{n\phi}{n-1} + sin^2\frac{n\phi}{n-1}}{\left(n^{\frac{n}{n-1}}\right)^2}\right)\\ - 2\left(\frac{cos\frac{\phi}{n-1}cos\frac{n\phi}{n-1} + sin\frac{\phi}{n-1}sin\frac{n\phi}{n-1}}{n^{\frac{1}{n-1}}n^{\frac{n}{n-1}}}\right)\nonumber
\end{aligned}
\end{equation}
\begin{equation}
\begin{aligned}
x^2 + y^2 = \left(\frac{1}{\left(n^{\frac{1}{n-1}}\right)^2}\right) + \left(\frac{1}{\left(n^{\frac{n}{n-1}}\right)^2}\right) \\- 2\left(\frac{cos\left(\frac{n\phi}{n-1} - \frac{\phi}{n-1}\right)}{n^{\frac{n+1}{n-1}}}\right)\nonumber
\end{aligned}
\end{equation}
\begin{equation}
\begin{aligned}
x^2 + y^2 = \left(\frac{1}{\left(n^{\frac{1}{n-1}}\right)^2}\right) + \left(\frac{1}{\left(n^{\frac{n}{n-1}}\right)^2}\right) - 2\left(\frac{cos\phi}{n^{\frac{n+1}{n-1}}}\right)
\end{aligned}
\label{eq.12}
\end{equation}

The above equation will have minima and maxima when $cos\Phi$ is $1$ and $-1$, respectively i.e. when $\Phi$ is $0$ and $\pi$. In order to correctly understand the shape of the corresponding Mandelbrot set of any polynomial of any degree, one needs to look at the above equation in conjunction with the corresponding parametric equations.

Since for $n=2$ the parametric equations have period of $2\pi$, the value of $x^2+y^2$ will have only one minimum i.e. at $\Phi=0$. For $n=3$, the period of corresponding parametric equations is $4\pi$ which means that to describe the full shape of the main body of the Mandelbrot set the parameter $\Phi$ will have to change from 0 to $4\pi$ during which $x^2 + y^2$ given by the above equation will go through two minima at $\Phi=0$ and $2\pi$ and two maxima at $\Phi$ values of $\pi$ and $3\pi$; for $n=4$, the period of the parametric equations is $6\pi$ which indicates three minima at $\Phi=0$, $2\pi$ and $4\pi$; and so on.

This argument can be safely extended to include that as n increases the number of minima (which are the indents) and maxima (which denote the mid-point of each cardioid) will increase and their number will be one less than the degree of the corresponding polynomial.

The plot of the central part of the Mandelbrot set for different values of $n$, as given in the Figure \ref{fig:card}, supports this point.

One more point that can be observed from these plots is that as n increases the size of the central part of the Mandelbrot set also gradually increases. At $n=100$, the radius of the central part is nearly $1$. Concept of limits can be used to prove that as n tends to infinity the expression,

\pagebreak
\begin{equation}
\begin{aligned}
x^2 + y^2 = \left(\frac{1}{\left(n^{\frac{1}{n-1}}\right)^2}\right) + \left(\frac{1}{\left(n^{\frac{n}{n-1}}\right)^2}\right) - \left(\frac{2cos\phi}{n^{\frac{n+1}{n-1}}}\right)
\end{aligned}
\label{eq13}
\end{equation}

will converge into an equation of a circle centered at the origin and with radius 1. As n increases and tends to infinity the second and the third term in Equation \ref{eq13}
vanish since the denominator contains $n^2$ and $n$ respectively together with $n^{\frac{2}{n-1}}$ (which is finite -- equals to $1$ -- as $n$ tends to infinity as will be shown below).

Now what remains is,

\begin{equation}
\begin{aligned}
x^2 + y^2 = \left(\frac{1}{(n^{\frac{1}{n-1}})^2}\right)
\end{aligned}
\end{equation}

It would be useful to examine the limit of $n^{\frac{1}{n-1}}$ as n tends to infinity. Using the fact that $A = e^{lnA}$, the expression $n^{\frac{1}{n-1}}$ , can be rewritten as $e^{ln(n^{\frac{1}{n-1}})}$. Now, since $e$ is a constant, the limit of $ln\left(n^{\frac{1}{n-1}}\right)$ as $n$ tends to infinity should be examined.

\begin{equation}
\begin{aligned}
\lim_{n \to \infty} ln\left(n^{\frac{1}{n-1}}\right) = \lim_{n \to \infty} \frac{ln(n)}{n-1}
\end{aligned}
\end{equation}

This limit assumes the form $\frac{\infty}{\infty}$ as $n$ tends to $\infty$. We can use L'hospital's rule to evaluate the limit. According to L'hospital's rule, $\lim_{x \to a} \frac{f(x)}{g(x)} = \lim_{x \to a} \frac{f'(x)}{g'(x)}$ if $\frac{f(x)}{g(x)}$ is either $\frac{0}{0}$ or $\frac{\infty}{\infty}$.

Hence, 
\begin{equation}
\begin{aligned}
\lim_{n \to \infty} ln\left(n^{\frac{1}{n-1}}\right) = \lim_{n \to \infty} \frac{ln(n)}{n-1} = \lim_{n \to \infty} \frac{\frac{1}{n}}{1} = 0\nonumber
\end{aligned}
\end{equation}
\begin{equation}
\begin{aligned}
\lim_{n \to \infty} n^{\frac{1}{n-1}} = 1
\end{aligned}
\end{equation}

In other words,

\begin{equation}
\begin{aligned}
\lim_{n \to \infty} \left(x^2 + y^2\right) = \left(\frac{1}{n^{\frac{1}{n-1}}}\right)^2 = 1^2 = 1
\end{aligned}
\end{equation}

Since $x^2 + y^2 = 1$ is a circle centered at the origin, shape of Mandelbrot set’s main body will tend to become a circle of unit radius as the power of the corresponding polynomial increases and tends to infinity. The changing shapes have already been seen in Figure \ref{fig:card} up to $n=100$.

\section{Convergence to a unit circle: Another approach}

Now that the limit of $\frac{1}{n^{\frac{1}{n-1}}}$ as $n$ tends to infinity has been deduced, another method can be discussed for determining the shape of a Mandelbrot set as the power of its corresponding polynomial function tends to infinity.

The square of the maximum distance, from the origin, of the points on the main body of the Mandelbrot sets of power $n$ can be derived from Equation \ref{eq.12}. As $cos\Phi$ can have values ranging from $1$ to $-1$, the maxima will occur when $cos\Phi$ is $-1$ and minima will occur when $cos\Phi$ is $1$.

\begin{equation}
\begin{aligned}
|c|_{max}^2 = \frac{1}{\left(n^{\frac{1}{n-1}}\right)^2} + \frac{1}{\left(n^{\frac{n}{n-1}}\right)^2} - \frac{2*(-1)}{n^{\frac{1}{n-1}}n^{\frac{n}{n-1}}}\nonumber
\end{aligned}
\end{equation}
\begin{equation}
\begin{aligned}
\implies |c|_{max}^2 = \frac{1}{\left(n^{\frac{1}{n-1}}\right)^2} + \frac{1}{\left(n^{\frac{n}{n-1}}\right)^2} + \frac{2}{n^{\frac{1}{n-1}}n^{\frac{n}{n-1}}}
\end{aligned}
\end{equation}

\begin{equation}
\begin{aligned}
|c|_{min}^2 = \frac{1}{\left(n^{\frac{1}{n-1}}\right)^2} + \frac{1}{\left(n^{\frac{n}{n-1}}\right)^2} - \frac{2*(1)}{n^{\frac{1}{n-1}}n^{\frac{n}{n-1}}}\nonumber
\end{aligned}
\end{equation}
\begin{equation}
\begin{aligned}
\implies |c|_{min}^2 = \frac{1}{\left(n^{\frac{1}{n-1}}\right)^2} + \frac{1}{\left(n^{\frac{n}{n-1}}\right)^2} - \frac{2}{n^{\frac{n+1}{n-1}}}
\end{aligned}
\end{equation}

Now, it would be useful to examine the value of the difference of the moduli of maxima and minima of the Mandelbrot set as $n$ tends to infinity.

\begin{equation}
\begin{aligned}
\lim_{n \to \infty} \left(|c|_{max}^2 - |c|_{min}^2\right) =\hfill \\\lim_{n \to \infty} \left(\frac{1}{\left(n^{\frac{1}{n-1}}\right)^2} + \frac{1}{\left(n^{\frac{n}{n-1}}\right)^2} + \frac{2}{n^{\frac{n+1}{n-1}}}\right) \\- \left(\frac{1}{\left(n^{\frac{1}{n-1}}\right)^2} + \frac{1}{\left(n^{\frac{n}{n-1}}\right)^2} - \frac{2}{n^{\frac{n+1}{n-1}}}\right)\nonumber
\end{aligned}
\end{equation}
\begin{equation}
\begin{aligned}
\implies \lim_{n \to \infty} \left(|c|_{max}^2 - |c|_{min}^2\right) = \lim_{n \to \infty} \frac{4}{n^{\frac{n+1}{n-1}}}\nonumber
\end{aligned}
\end{equation}
\begin{equation}
\begin{aligned}
\implies \lim_{n \to \infty} \left(|c|_{max}^2 - |c|_{min}^2\right) = \lim_{n \to \infty} \left(\frac{1}{n^{\frac{1}{n-1}}}\right)^2\left(\frac{4}{n}\right)\nonumber
\end{aligned}
\end{equation}
\begin{equation}
\begin{aligned}
\implies \lim_{n \to \infty} \left(|c|_{max}^2 - |c|_{min}^2\right) = \lim_{n \to \infty} \left(\frac{1}{n^{\frac{1}{n-1}}}\right)^2\lim_{n \to \infty}\left(\frac{4}{n}\right)\nonumber
\end{aligned}
\end{equation}
\begin{equation}
\begin{aligned}
\implies \lim_{n \to \infty} \left(|c|_{max}^2 - |c|_{min}^2\right) = \lim_{n \to \infty}\left(\frac{1}{1}\right)^2\left(\frac{4}{n}\right)\nonumber
\end{aligned}
\end{equation}
\begin{equation}
\begin{aligned}
\implies \lim_{n \to \infty} \left(|c|_{max}^2 - |c|_{min}^2\right) = 0\nonumber
\end{aligned}
\end{equation}
\begin{equation}
\begin{aligned}
\implies \lim_{n \to \infty} \left(|c|_{max} - |c|_{min}\right) = 0
\end{aligned}
\end{equation}

Since the limit of $|c|_{max} - |c|_{min}$ is zero as $n$ tends to infinity, it can be deduced that the maxima and minima of c values tend to be equidistant from the origin as the power of its corresponding polynomial function tends to infinity. Hence, the shape of a Mandelbrot set tends towards a circle as the power of its polynomial function tends to infinity.

\section{Connection with roots of unity}

It is observed here that as $n$ is increased, the number of cardioids also increases. For $n=2$ there is only one cardioid which is the central part of the Mandelbrot set. For $n=3$ there are two cardioids connected back-to-back. Similarly, for $n=4$ and $n=5$ there are three and four cardioids respectively connected back-to-back. For the general case, for an $n$ degree iterated function that there would be $(n-1)$ cardioids connected back-to-back.

This conclusion is closely related to the fact that for any iterated function the fixed points are determined by the differential criterion. Since the first differentiation of any polynomial function reduces its degree by one, the number of fixed points determined by the differentiation criterion will be one less than the degree of function.

For example, for the following function,
\begin{equation}
\begin{aligned}
f(z) = z^n + c\nonumber
\end{aligned}
\end{equation}
\begin{equation}
\begin{aligned}
f'(z) = nz^{n-1}
\end{aligned}
\end{equation}

The differentiation criterion yields $nz^{n-1} = 1$. This gives $z = (\frac{1}{n})^{\frac{1}{n-1}}cis\frac{\Phi}{n-1}$. This value of $z$ satisfying the differentiation criterion is then substituted in the equation $z = z^n + c$ to get the positions of the complex number $c \left(= x + iy\right)$ in terms of $\Phi$.

In order to determine the location of maximum/minimum differential calculus may be used. Function $z = z^n + c$ may be written as $c = z – z^n$. According to differential criterion for stationary points $\left(maxima/minima\right)$ the first order differentiation is equal to zero. Differentiating $c$ with respect to $z$:

\begin{equation}
\begin{aligned}
\frac{dc}{dz} = 1 - nz^{n-1}
\end{aligned}
\end{equation}

If $\frac{dc}{dz}$ is equal to zero, $1 - nz^{n-1} = 0$, which is the differential criterion of fixed points. In other words, the points of the plot of parameter $c$ where it has stationary values will be given by the roots of the equation $nz^{n-1} = 1$.

Substituting this value of $z^{n-1}$ in the following equation,

\begin{equation}
\begin{aligned}
c = z - z^n\nonumber
\end{aligned}
\end{equation}
\begin{equation}
\begin{aligned}
c = z\left(1 - z^{n-1}\right)\nonumber
\end{aligned}
\end{equation}
\begin{equation}
\begin{aligned}
c = z\left(1 - \frac{1}{n}\right)\nonumber
\end{aligned}
\end{equation}

It is obvious from above that argument of minimum values of $c$ is same as that of $z$.

The value of $z$ from the differential criterion of fixed points is $z = \left(\frac{1}{n}\right)^{\frac{1}{n-1}}cis\frac{\Phi}{n-1}$. Geometrically this equation means that argument of $z$ has $\left(n-1\right)$ values which are equally spaced on the unit circle in the Argand plane. In other words, the roots of unity will provide the arguments of the points where the plots of the parameter $c$, which form the main body of the Mandelbrot set, will have stationary values.

Applying this argument to the biquadratic iterated function $z^4 + c$ one gets $4z^3 = 1$ or $z = \left(\frac{1}{4}\right)^{\frac{1}{3}}$, $\left(\frac{1}{4}\right)^{\frac{1}{3}}\omega$, $\left(\frac{1}{4}\right)^{\frac{1}{3}}\omega^2$ as three different values of $z$ where $\omega$ denotes one of the cube roots of unity; the others being $1$ and $\omega^2$.

Substituting $z^3 = \frac{1}{4}$ in the following equation,

\begin{equation}
\begin{aligned}
c = z – z^4 = z\left(1 - z^3\right)\nonumber
\end{aligned}
\end{equation}
\begin{equation}
\begin{aligned}
c = z\left(1 - \frac{1}{4}\right) = \frac{3z}{4}
\end{aligned}
\end{equation}

Thus, minimum values of $c$ will have arguments of $1$, $\omega$, and $\omega^2$ which are the cube roots of unity. This conclusion is also supported by the experimental plot of the Mandelbrot set which has three Cardioids joined back-to-back. Hence, one has a simple yet elegant way to locate the indent points of the main body of the Mandelbrot set corresponding to a polynomial of any positive integer power

\section{Conclusion}

From the preceding paragraphs, the direct relationship between the geometric shape of the central part of the Mandelbrot set and the power of the corresponding polynomial function has become clear. It is observed that for an iterated complex function of positive integral power the number of cardioids in the central part of the corresponding Mandelbrot is one less than the degree of the polynomial. The property of the points located on the boundary of the central part, that their derivative equals $1$, makes it possible to derive parametric equations for the geometric shape. As the positive integral power of the polynomial function tends to infinity, the central part of Mandelbrot sets is proved to tend towards a unit circle. It is also seen that the arguments of the points where the cardioids meet are indicated by the roots of the unity. This analysis can be extended to cover the cases of fractional power of the iterated functions.

\section*{Acknowledgement}
I would like to thank my mathematics teacher Dovid Fein for his supervision and support.


\begin{thebibliography}{9}
  
\bibitem{brown2003}
Brown, James Ward, and Ruel Vance Churchill. \textit{"Complex variables and applications."} McGraw-Hill. 2003.
  
\bibitem{carleson2005}
Carleson, Lennart, and Theodore W. Gamelin. \textit{"Complex dynamics."} Springer-Verlag. 2005.

\bibitem{dickerson2003}
Dickerson, Richard. \textit{"Higher Order Mandelbrot Fractals: Experiments in Nanogeometry."} 2003.

\bibitem{falconer1990}
Falconer, Kenneth, John Wiley. \textit{"Fractal Geometry: Mathematical Foundations and Applications"} 1990.

\bibitem{gamelin2000}
Gamelin, Theodore W. \textit{"Complex Analysis."} Springer. 2000.

\bibitem{mathews2012}
Mathews, John H., and Russell W. Howell. \textit{"Complex analysis for mathematics and engineering."} Jones
and Bartlett, 2012.

\bibitem{peitgen1992}
Peitgen, Heinz-Otto, Hartmut Jurgens and Dietmar Saupe. \textit{"Chaos and Fractals: New Frontiers of Science."} 1992.

\bibitem{stewart1999}
Stewart, James. \textit{"Calculus."} Gary W. Ostedt. 1999.

\end{thebibliography}
\end{document}